\documentclass[a4paper,english,12pt]{article}
\usepackage[logonly]{trace}
\usepackage{babel}
\usepackage{amsfonts, amsmath, amssymb}
\usepackage{graphicx}
\usepackage{pstricks}
\usepackage{psfig}
\usepackage{multirow}
\usepackage{fancyhdr}
\usepackage{vmargin, fancybox}

\newcommand{\mathsym}[1]{{}}

\begin{document}

\title{A Note On ``Solitary Wave Solutions of the Compound Burgers-Korteweg-de
Vries Equation"}

\author{Claire David\footnotemark[2] \footnote{Corresponding author:
david@lmm.jussieu.fr; fax number: (+33) 1.44.27.52.59.} ,  Rasika
Fernando\footnotemark[2] , and Zhaosheng Feng\footnotemark[3]
\\ \\
\footnotemark[2] \small{Universit\'e Pierre et Marie Curie-Paris 6}  \\
\small{Laboratoire de Mod\'elisation en M\'ecanique, UMR CNRS 7607,} \\
\small{Bo\^ite courrier $n^0162$, 4 place
Jussieu, 75252 Paris, cedex 05, France} \\
\footnotemark[3] \small{Department of Mathematics, University of
Texas-Pan American, Edinburg, TX 78541, USA} \\}

\date{}

\maketitle

\begin{abstract}
The goal of this note is to construct a class of traveling
solitary wave solutions for the compound Burgers-Korteweg-de Vries
equation by means of a hyperbolic ansatz. A computational error in
a previous work has been clarified.
\end{abstract}

\section{Introduction}
\label{sec:intro}

\noindent

Consider the following equations:\\ \\
$\diamond$ KdV:
\begin{equation} u_t + \alpha u u_x + s u_{xxx} = 0, \end{equation}
$\diamond$ Burgers: \begin{equation}  u_t + \alpha u u_x + \mu
u_{xx} = 0, \end{equation} $\diamond$ modified KdV:
\begin{equation}  u_t + \beta u^2 u_x + s u_{xxx} = 0 \, \,\,(s>0) \end{equation}
These three equations play crucial roles in the history of wave
equations. Equation (1) is named after its use by Burgers
\cite{burger1} for studying turbulence in 1939. Equation (2) was
first suggested by Korteweg and de Vries \cite{kdv1} who used it as
a nonlinear model to study the change of form of long waves
advancing in a rectangular channel. It was shown by Wadati
\cite{wadati} that equation (3) can be solved exactly by the inverse
scattering method and the N-soliton solution can be expressed
explicitly which is essentially the same as that in the KdV
equation. A combination of these three equations leads to the
so-called compound Burgers-Korteweg-de Vries (cBKDV) equation:
\begin{equation}
\label{CBKDV} u_t + \alpha u u_x + \beta u^2 u_x + \mu u_{xx} + s
u_{xxx} = 0,
\end{equation}
where $\alpha$, $\beta$, $\mu$ and $s$ are real constants. The first
term is the \emph{instationary term}. The second and third ones are
two \emph{non-linear convective terms} with different orders. The
fourth is the so-called \emph{viscous dissipative term}: $\mu $,
which refers to a positive quantity, denotes the viscosity. The last
one is the \emph{dispersive term}. In order for the traveling wave
to keep its appearance all through its propagation, an equilibrium
is required between convective and dispersive terms.
\\
\noindent The coefficient $\alpha$, which usually refers to the wave
celerity, will be assumed positive in the following. The case
$\alpha<0$ can be easily deduced from the latter one by
adimensionalization.\\

\noindent The coefficient of the $u_{xxx}$ term in the modified KDV
equation being strictly positive, we restrain our study to the case
$s>0$ in (\ref{CBKDV}). As in [8], the quantity $-s\,\beta$ is
supposed positive, which results in $\beta<0$.\\

In the past few years, the cBKDV equation has attracted continuous
attention from a diverse group of researchers. Especially, the
traveling wave solution has been explored extensively. An exact
solution was presented by Wang \cite{wang1} using the homogeneous
balance method. Feng [5-7] investigated the same problem by
applying the first integral method and the method of variation of
parameters, respectively. A more general result was obtained by
Parkes and Duffy \cite{park1, park2} using the automated
tanh-function method and adapting the method of variation of
parameters used in \cite{feng3}. A generalization of equation (4)
with higher-order nonlinearities in one- and two- dimensional
spaces was treated by Zhang et al. \cite{zhang1, zhang2} using the
method of undetermined coefficients, and by Li et al. \cite{li}
applying the tanh-function method and direct assumption method
with symbolic computations, respectively. Describing traveling
waves of nonlinear evolution equations has been one of basic
problems in theoretical and experimental physics. Traveling wave
solutions to many nonlinear one-dimensional evolution equations
can be derived from a set of ordinary differential equations that
can be interpreted as a flow in a phase space. In the last
century, experiments on one-dimensional states of nonlinear
traveling wave convection were undertook by many researchers in
the narrow annular cell, condensed matter physics, plasma,
particle physics, and so on [13-16].
\\

In the present work, we aim at extending results in a previous
work \cite{feng1}, in which traveling wave solutions of the cBKDV
equation were exhibited as combinations of bell-profile waves and
kink-profile waves. Following [5], we assume that equation (4) has
the traveling wave solution of the form
\begin{equation} u(x, t) =u(\xi), \quad \xi=x-vt \end{equation}
where $v$ is the wave velocity. Substituting(5) into equation (4),
performing an integration with respect to $\xi$ and setting the
integration constant to zero yields
\begin{equation}
u''( \xi ) + ru'( \xi ) + au( \xi )^3 + bu( \xi )^2 + cu( \xi ) = 0,
\end{equation}
where $r=\frac{\mu}{s}$, $a=\frac{\beta}{3s}$,
$b=\frac{\alpha}{2s}$ and $c=-\frac{v}{s}$. Equation (6) is a
second-order nonlinear ordinary differential equation. We know
that if all coefficients of equation (4) are nonzero, equation (6)
does not pass the Painlev\'e test [17, 18].
\\

The rest of the paper is organized as follows. In Section 2, a
class of traveling wave solutions to the cBKDV equation is
presented by using a hyperbolic ansatz and the types of traveling
waves are discussed theoretically. Numerical simulations for a
couple of sets of parameters are illustrated. A calculational
error in the previous literature is clarified at the end of this
section. In Section 3, we give a brief discussion.

\section{Traveling Solitary Waves}

\subsection{Hyperbolic Ansatz}

\noindent The phase plane analysis in Section 3 of Paper [5]
provides us useful information when we construct traveling solitary
wave solutions for equation (4). It clearly indicates that under
certain parametric conditions, equation (4) does have bounded
kink-profile traveling wave solutions. Based on this result, in this
section, a class of traveling wave solutions of the cBKDV equation
is established by using a hyperbolic ansatz which is actually a
combination of bell-profile waves and kink-profile waves of the form
\begin{equation}
u(x, t) = \sum_{i = 1}^n \left (B_i\; \text{tanh}^i \left [C_i
(x-vt+x_0) \right ] + D_i \; \text{sech}^i \left [C_i(x-vt+x_0)
\right ] \right ) + B_0,
\end{equation}
where the $B_i's$, $C_i's$ $D_i's$, $(i=1,\ \cdots,\ n)$, $v$ and
$B_0$ are constants to be determined, and $x_0$ is arbitrary.
\\

\noindent After substitution of (7) into equation (6), we get
\begin{equation}
\begin{array}{l}
\scriptsize{ \underset{i=1}{\overset{n}{\sum }} \left \lbrace
\begin{array}{c}
i \;(-\text{sech}[C_i\,\xi]^{i} +(1+i)\; \text{sech}[
C_i\,\xi]^{2+i}\; \text{sinh}[z C_i]^2) C_i^{2} D_i\\+ B_i C_i^{2}
((-1+i) i\; \text{sech}[C_i\,\xi]^4 \;\text{tanh}[C_i\,\xi]^{-2+i}
-2 i \;\text{sech}[C_i\,\xi]^2 \;\text{tanh}[ C_i\,\xi]^i)
\end{array}
\right \rbrace} \\
   \\
+\frac{\mu \; \underset{i=1}{\overset{n}{\sum }}(-i\;
\text{sech}[C_i\,\xi]^{1+i} \text{sinh}[C_i\,\xi] C_i D_{i}+i\;
\text{Sech}[C_i\,\xi]^2 B_i C_i \text{tanh}[z C_i]^{-1+i})}
{s}\\
+\frac{\beta \;(\underset{i = 1}{\overset{n}{\sum }}
(B_i\;\text{tanh}^i[C_i\,\xi]+D_i\;
\text{sech}^i[C_i\,\xi])+B_0)^3}{3s}\\
+\frac{\alpha \;(\underset{i = 1}{\overset{n}{\sum }}
(B_i\;\text{tanh}^i[C_i\,\xi]+D_i\;\text{Sech}^i[C_i\,\xi])+B_0)^2}{2s}
-\frac{v \;(\underset{i = 1}{\overset{n}{\sum }}
(B_i\;\text{tanh}^i[C_i\,\xi]+D_i\;\text{sech}^i[C_i\,\xi])+B_0)}{s}=0.
\end{array}
\end{equation}
\noindent The difficulty for solving this equation lies in finding
the values of the constants $B_i$, $C_i$,  $D_i$, $B_0$ and $v$ by
using the over-determined algebraic equations. Following
\cite{feng1}, after balancing the higher-order derivative term and
the leading nonlinear term, we deduce $n=1$. Then we replace
$\mbox{sech}(C_1\,\xi)$ by $ \frac{2}{  e^{ {C_1}\,\xi}+e^{-
{C_1}\,\xi } }$, $\mbox{sinh(}C_1\,\xi)$ by $ \frac {
e^{{C_1}\,\xi}-e^{- {C_1}\,\xi } }{2}$, $\mbox{tanh}(C_1\,\xi)$ by
$\frac{e^{  {C_1}\,\xi
 } -
      e^{-  {C_1}\,\xi}}{e^{  {C_1}\,\xi  } + e^{-
      {C_1}\,\xi}}$, and multiply both sides by ${({  e^{ {C_1}\,\xi}+e^{-
{C_1}\,\xi  } })}\, e^{ {C_1}\,\xi}$, so equation (8) can be
rewritten in the following form:
\begin {equation}
\sum_{k=0}^{6} P_k \left ( B_0,\ B_1,\ C_1,\ D_1,\ v \right) e^{k\,
C_1 \,\xi} \; = \;0,
\end{equation}
where the $P_k$ $(k=0,\,...,\,6 \,)$, are polynomials of $B_0$,
$B_1$, $C_1$, $D_1$ and $v$.
\\

\noindent With the aid of mathematical softwares such as
Mathematica, when $n=1$, equating the coefficient of each term $e^{k
\,C_1\,\xi}$ $(k=0, ..., 6\,)$ in equation (9) to zero, yields a
nonlinear algebraic system which contains seven equations as
follows:
\begin {equation}
\scriptsize{\begin{array}{rcl}

   P_0 & = &  \scriptsize {\left  ( \begin{array}{rcl}

   &  &

   -24\,
  v\, {B_0} + 12\, \alpha \, {{B_0}}^2 - 8\, \beta \, {{B_0}}^3 +
    24\, \mu \, {B_1}\, {C_1} - 24\,
  v\, {D_1} +24\, \alpha \, {B_0}\, {D_1} +
    24\, \beta \, {{B_0}}^2\, {D_1}  \\

    & - &  24\,
  s\, {{C_1}}^2\, {D_1} + 12\, \alpha \, {{D_1}}^2 +
    24\, \beta \, {B_0}\, {{D_1}}^2 + 8\, \beta \, {{D_1}}^3

    \end{array}   \right )=0,}
    \\

 &   &     \\

  P_1 & = &  \scriptsize {  \left  ( \begin{array}{rcl}

     &  &  -6\, v\, {D_1} + 6\, \alpha \, {B_0}\, {D_1} +
  6\, \beta \, {{B_0}}^2\, {D_1} - 6\, \alpha \, {B_1}\, {D_1} -
  12\, \beta \, {B_0}\, {B_1}\, {D_1}  \\

  & + &   6\, \beta \, {{B_1}}^2\, {D_1} +
  6\, \mu \, {C_1}\, {D_1} + 6\, s\, {{C_1}}^2\, {D_1}

     \end{array}

    \right )=0,  }  \\

     &   &     \\

   P_2 &  = &  \scriptsize{  \left  ( \begin{array}{rcl}

 &   &   -9\,
  v\, {B_0} + \frac{9\, \alpha \, {{B_0}}^2}{2} + 3\, \beta \, {{B_0}}^3
  +
    3\, v\, {B_1} - 3\, \alpha \, {B_0}\, {B_1} -
    3\, \beta \, {{B_0}}^2\, {B_1} - \frac{3\, \alpha \,
    {{B_1}}^2}{2}  \\

   & - &
    3\, \beta \, {B_0}\, {{B_1}}^2 + 3\, \beta \, {{B_1}}^3 +
    12\, \mu \, {B_1}\, {C_1} + 24\,
  s\, {B_1}\, {{C_1}}^2 + 6\, \alpha \, {{D_1}}^2  \\

  & + &
    12\, \beta \, {B_0}\, {{D_1}}^2 - 12\, \beta \, {B_1}\, {{D_1}}^2  \end{array}

    \right )=0, }\\

 &   &     \\

    P_3 & = &  \scriptsize{   -12\,
  v\, {D_1} + 12\, \alpha \, {B_0}\, {D_1} +
    12\, \beta \, {{B_0}}^2\, {D_1} - 12\, \beta \, {{B_1}}^2\, {D_1} - 36\,
  s\, {{C_1}}^2\, {D_1} + 8\, \beta \, {{D_1}}^3 }=0,
      \\
%    & =& 0,\\

 &   &     \\

   P_4 & = &  \scriptsize{  \left  ( \begin{array}{rcl}

     & - & 9\,
  v\, {B_0} + \frac{9\, \alpha \, {{B_0}}^2}{2} + 3\, \beta \, {{B_0}}^3 -
    3\, v\, {B_1} + 3\, \alpha \, {B_0}\, {B_1} +
    3\, \beta \, {{B_0}}^2\, {B_1} - \frac{3\, \alpha \,
    {{B_1}}^2}{2}  \\

   & - &
    3\, \beta \, {B_0}\, {{B_1}}^2 - 3\, \beta \, {{B_1}}^3 +
    12\, \mu \, {B_1}\, {C_1} - 24\,
  s\, {B_1}\, {{C_1}}^2 + 6\, \alpha \, {{D_1}}^2 +
    12\, \beta \, {B_0}\, {{D_1}}^2  \\
     & +

 & 12\, \beta \, {B_1}\,
    {{D_1}}^2 \end{array}

    \right )=0, }

   \\
    &   &     \\

   P_5 & = & \scriptsize{   \left  ( \begin{array}{rcl}

    &  &  -6\, v\, {D_1} + 6\, \alpha \, {B_0}\, {D_1} +
  6\, \beta \, {{B_0}}^2\, {D_1} + 6\, \alpha \, {B_1}\, {D_1} +
  12\, \beta \, {B_0}\, {B_1}\, {D_1}  \\

  &  + &   6\, \beta \, {{B_1}}^2\, {D_1} -
  6\, \mu \, {C_1}\, {D_1} + 6\, s\, {{C_1}}^2\, {D_1} \end{array}

    \right )=0, }\scriptsize

   \\
    &   &     \\

   P_6 & = &  \scriptsize{  \left  ( \begin{array}{rcl}
  & -&  3\,
  v\, {B_0} + \frac{3\, \alpha \, {{B_0}}^2}{2} + \beta \, {{B_0}}^3 - 3\,
  v\, {B_1} + 3\, \alpha \, {B_0}\, {B_1} +
    3\, \beta \, {{B_0}}^2\, {B_1} + \frac{3\, \alpha \, {{B_1}}^2}{2}

   \\
  & + &
    3\, \beta \, {B_0}\, {{B_1}}^2 + \beta \, {{B_1}}^3 \end{array}

    \right )=0.  }\scriptsize
\end{array}}
\end{equation}
\noindent System (10) can be solved consistently by using
Mathematica again.\\
\noindent For sake of simplicity, we use $\varepsilon_1$,
$\varepsilon_2$, $\varepsilon_3$  and $\varepsilon$ to denote $1$ or
$-1$, and denote by $\kappa$ the quantity:

\begin{equation}
\kappa=\varepsilon_1\,\frac{\alpha}{2\,|\beta|}\,\sqrt{\frac{|\beta|}{6\,
s }}-\varepsilon_2\,\frac{\mu}{6\, s }
\end{equation}

\noindent The sets of solutions are given by:

\begin{equation}
\left \lbrace
\begin{array}{rcl}
B_{0} &  =   & -\frac{\alpha}{2\,\beta}- \varepsilon_3\,
\;\sqrt{\frac{6\,s}{|\beta|}}\,\frac{\mu}{6\,s}
\\
B_{1}  &  =   & \varepsilon_3\, \sqrt{\frac{6\,s}{|\beta|}}\,\kappa \\
C_{1}  &  =   & 2\,\kappa
  \\
D_{1}  &  =   &   i\,\varepsilon \big
     (\frac{\alpha}{2\,|\beta|}-\varepsilon_3\,\frac {\mu}{6\,s}\,\sqrt{\frac{6\,s}{|\beta|}}  \,\big )\\
     v &  =   & - \frac{\mu^2}{6\,s}-2\,s\,\big
     [\frac{\alpha}{2\,|\beta|}\,\sqrt{\frac
     {|\beta|}{6\,s}}-\varepsilon_3\,\frac {\mu}{6\,s} \big ]^2-  \frac{\alpha^2
     }{4\,\beta}
\end{array}
\right. \end{equation}

\noindent with the constraint:
\begin{equation}
\label{contraint} \varepsilon_1\,\varepsilon_2\,\varepsilon_3=1
 \end{equation}
 \noindent or:

\begin{equation}
B_1=D_1=0
\end{equation}

\noindent Note that the second set to the fourth set of solutions
are not of much interest, since they only correspond to three
trivial cases---the constant solution. Now we examine the first set
of solutions. It is interesting to notice that the nature of
traveling wave solutions completely depends on the values of the
parameters $\alpha$, $\beta$, $s$, $\mu$. Moreover, we observe that
\begin{enumerate}
\item[\it{i.}] depending on the sign of $\frac{\mu^2}{6\,s}+2\,s\,\big
     [\frac{\alpha}{2\,\beta}\,\sqrt{\frac
     {|\beta|}{6\,s}}-\varepsilon_3\,\frac {\mu}{6\,s} \big ]^2-  \frac{\alpha^2
     }{4\,\beta}$, the traveling
wave solution moves positively or negatively;

\item[\it{ii.}] the case of $B_1 \neq 0$ corresponds
to a kink-profile wave solution;

\item[\it{iii.}] since $-s\,\beta>0$, $D_1$ is purely imaginary; equation (4) thus admits a complex
traveling solitary wave solution.\\
\noindent Due to the $tanh $ function, the real part of the solution
$u$ has a kink-profile, while its imaginary part has a bell-profile
(due to the $sech $ function).

\end{enumerate}

\subsection{Parametric Study}

\noindent In the following, we are going to show the influence of
parameters $\alpha$, $\beta$, $s$, $\mu$ on the shape and velocity
of traveling waves, and a couple of traveling waves corresponding to
different sets of parameters are illustrated.

\subsubsection{Shape of Traveling Waves}

\noindent The variations of the coefficients $B_1$ and $D_1$ have
crucial influence on the bell-profile and kink-profile of traveling
solitary
waves.\\

\noindent For this purpose, we consider the quotient of $B_1^2$ and
$D_1^2$:
\begin{equation}
\frac{{B_1}^2}{{D_1}^2}=
-\frac{\big(\frac{\alpha}{2\,|\beta|}-\varepsilon_1\,\varepsilon_2\,\frac{\mu}{\sqrt{\,6\,s\,|\beta|}}
\big)^2}{\big(\frac{\alpha}{2\,|\beta|}-\varepsilon_3\,\frac{\mu}{\sqrt{\,6\,s\,|\beta|}}
\big)^2}
\end{equation}

\noindent Due to the contraint (\ref{contraint}):

\begin{equation}
\varepsilon _3=\varepsilon _1\,\varepsilon _2
\end{equation}

\noindent It ensures:

\begin{equation}
{B_1}^2=-{D_1}^2
\end{equation}

\noindent and, since $D_1$ is pure imaginary:

\begin{equation}|B_1|=|D_1|
\end{equation}

\noindent which results in an equilibrium between the bell-profile
and kink-profile of the solitary wave.

\subsubsection{Velocity Analysis}

\noindent Here we extend our attention to the wave velocity $v$.
From the formula of $v$ in the first set of solutions, we can see
the rate of change of $v$ (with negative option) with respect to
$\alpha$ and $\mu$, respectively:
\begin{equation*}
\frac{\partial {v}} {\partial \alpha}=\frac{
\,\sqrt{\frac{6\,|\beta| }{s}}\, \varepsilon_3\, \mu +6 \, \alpha
}{18\, |\beta| }
\end{equation*}
\begin{equation*}
\frac{\partial {v}} {\partial \mu}=\frac{ \,s \,\alpha\,
\sqrt{\frac{6\,\beta }{s}} \,\varepsilon_3 -8\, |\beta|
   \, \mu }{18 \,s\, |\beta| }
\end{equation*}
\noindent Since $\beta<0$, there exists critical points for only one
of those derivatives. For example, if $\varepsilon_3=-1$, denote:
\begin{equation}
 \alpha_v=\mu\,{\sqrt{\frac{|\beta| }{6\,s}}  }
\,\,\, , \,\,\,
\end{equation}
\noindent In this specific case, we can derive variational tables as
follows directly
\begin{equation*}
\begin{array}{|c|ccccr|}
\hline
\alpha     & 0  &    &  {\alpha _v}  &     & +\infty  \\
\hline
\frac{\partial {v}} {\partial \alpha } &          &  - & 0  & + &   \\
\hline
 & -\frac{2\,\mu ^2}{9 \,s } & & & & +\infty \\       % ligne des valeurs "max"
v &  &\searrow   &  &\nearrow  & \\   % flèches
&  & & -\frac{\mu ^2}{4\, s}    & &  \\          % ligne des valeurs "min"
\hline
\end{array}\,\,\,  \,\,\,
\begin{array}{|c|cccr|}
\hline
\mu    & 0  &    &     & +\infty  \\
\hline
\frac{\partial {v}} {\partial \mu } &          &  -   & &   \\
\hline
 & \frac{\alpha ^2}{6 \,|\beta| } &  & &  \\       % ligne des valeurs "max"
v &  &\searrow   &    & \\   % flèches
&  &      & & -\infty \\          % ligne des valeurs "min"
\hline
\end{array}
\end{equation*}
%\noindent For $\alpha=\alpha_c=\frac{\sqrt{\frac{2}{3}\, \beta }
%\mu }{\sqrt{s}}$ or
%$\mu=\mu_c=\frac{\sqrt{\frac{3}{2}\,s} \alpha }{\sqrt{\beta }}$, $v$ vanishes.\\

\noindent Denote respectively $\alpha_c$ the value of the parameter
$\alpha$ which satisfies $\alpha_v<\alpha_c$, $v(\alpha_c)=0$, and
$\mu_c$ the value of the parameter $\mu$ which satisfies $0<\mu_c$,
$v(\mu_c)=0$. When $\alpha$ varies in $[\alpha_c,+\infty[$ or $\mu$
varies in $[0,\mu_c]$, the wave will propagate with a positive
velocity. When $\alpha$ gets larger and larger alone, the wave will
propagate with a big positive speed; but when $\mu$ becomes larger
and larger alone, the wave will propagate with a big negative speed.
\\

\noindent Similarly, from the derivative of $v$ with respect to
$\beta$ and $s$
\begin{equation*}
\frac{\partial {v}} {\partial |\beta|}=-\frac{\alpha \, \left(6
\,\alpha +\sqrt{\frac{6\,|\beta| }{s}}\, \varepsilon_3\,
   \mu \right)}{36\, \beta ^2} \,\,\, ,  \,\,\,
   \frac{\partial {v}} {\partial s}=\frac{\mu  \left(-\sqrt{6}\, \alpha \, \varepsilon_3+8\,\mu\, \sqrt{\frac{|\beta| }{s}}
   \right)}{36\, s^2 \,\sqrt{\frac{| \beta|  }{s}}}
\end{equation*}

\noindent There exists critical points for only one of those
derivatives. For example, if $\varepsilon_3=-1$, denote:
\begin{equation}
\beta_v=-\frac{\alpha^2\,s }{6\,\mu^2 } \end{equation} We can find
the following table immediately
%\noindent $\frac{\partial {v}} {\partial \beta}$ never vanish, $
%\frac{\partial {v}} {\partial \mu}$
%vanishes for $s={s_v}=\frac{32 \beta  \mu ^2}{3 \alpha ^2}$.\\
%\\
%\noindent The following variation tables are thus obtained:
\begin{equation*}
\begin{array}{|c|ccccr|}
\hline
|\beta |  & 0   &    &  |{\beta_v}|  &     & +\infty \\
\hline
\frac{\partial {v}} {\partial | \beta |} &          &  - & 0  & + &   \\
\hline
 & +\infty & &  & & -\frac{\mu ^2}{9 \,s }\\       % ligne des valeurs "max"
v &  &\searrow   &  &\nearrow  & \\   % flèches
&  & &  -\scriptsize{\frac{\mu^2}{4 \,s }}  & &  \\          % ligne des valeurs "min"
\hline
\end{array}
\,\,\,  \,\,\,
\begin{array}{|c|cccr|}
\hline
s   & 0  &    &     & +\infty   \\
\hline
\frac{\partial {v}} {\partial s} &          &  -   &  &   \\
\hline
 &  & &  &-\frac{\alpha ^2}{4 \,\beta }\\       % ligne des valeurs "max"
v &     &  &\nearrow  & \\   % flèches
&
  -\infty  & & &  \\          % ligne des valeurs "min"
\hline
\end{array}
\end{equation*}
\noindent As $\beta$ approaches to zero, the absolute value of the
wave speed will increase and eventually blows up. Similar thing will
occur when $s$ approaches to zero, but the wave will propagate with
a positive speed.

\subsubsection{Numerical Example}

\begin{figure}[h!]
\centerline{
 \includegraphics[width=6cm]{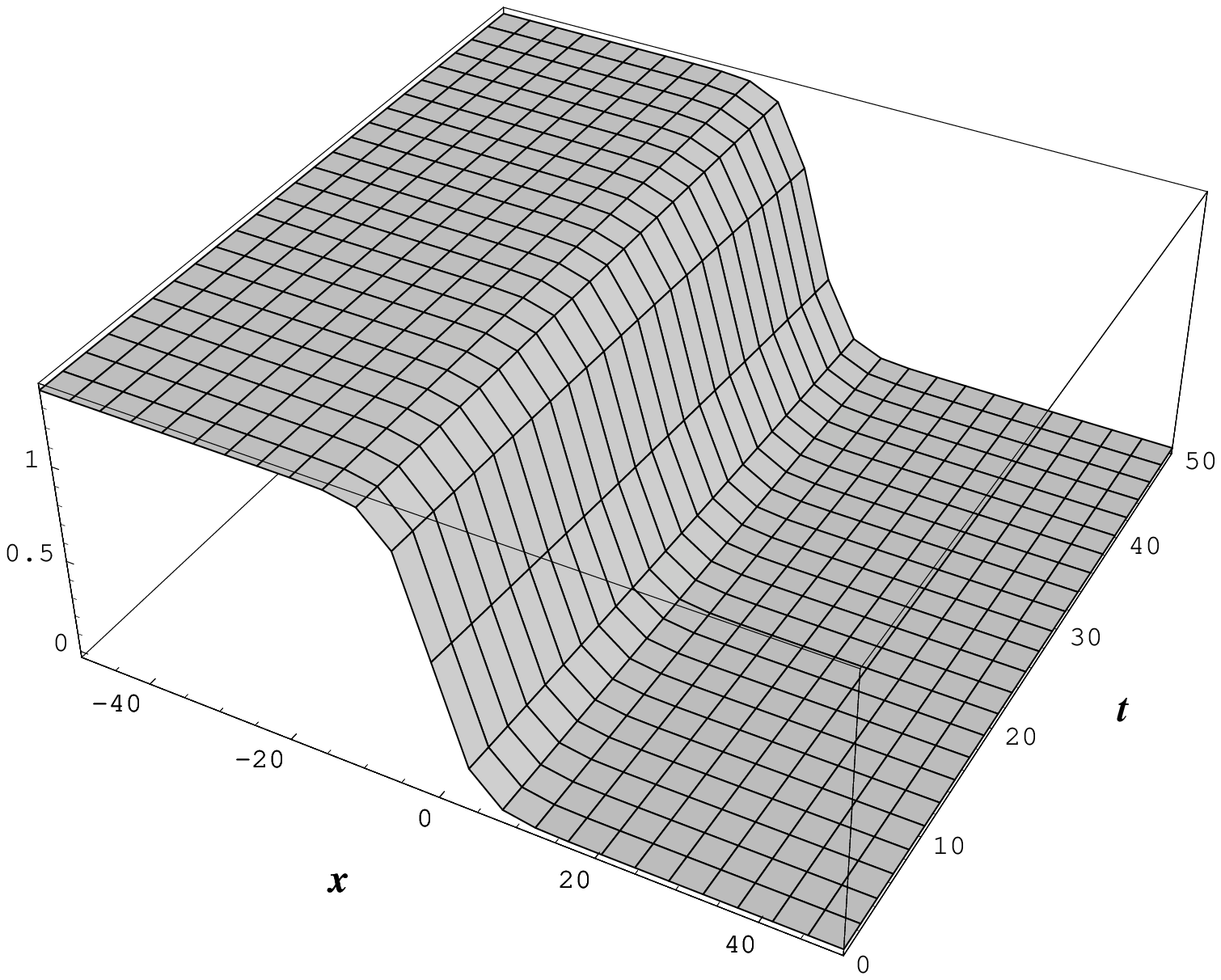} \includegraphics[width=6cm]{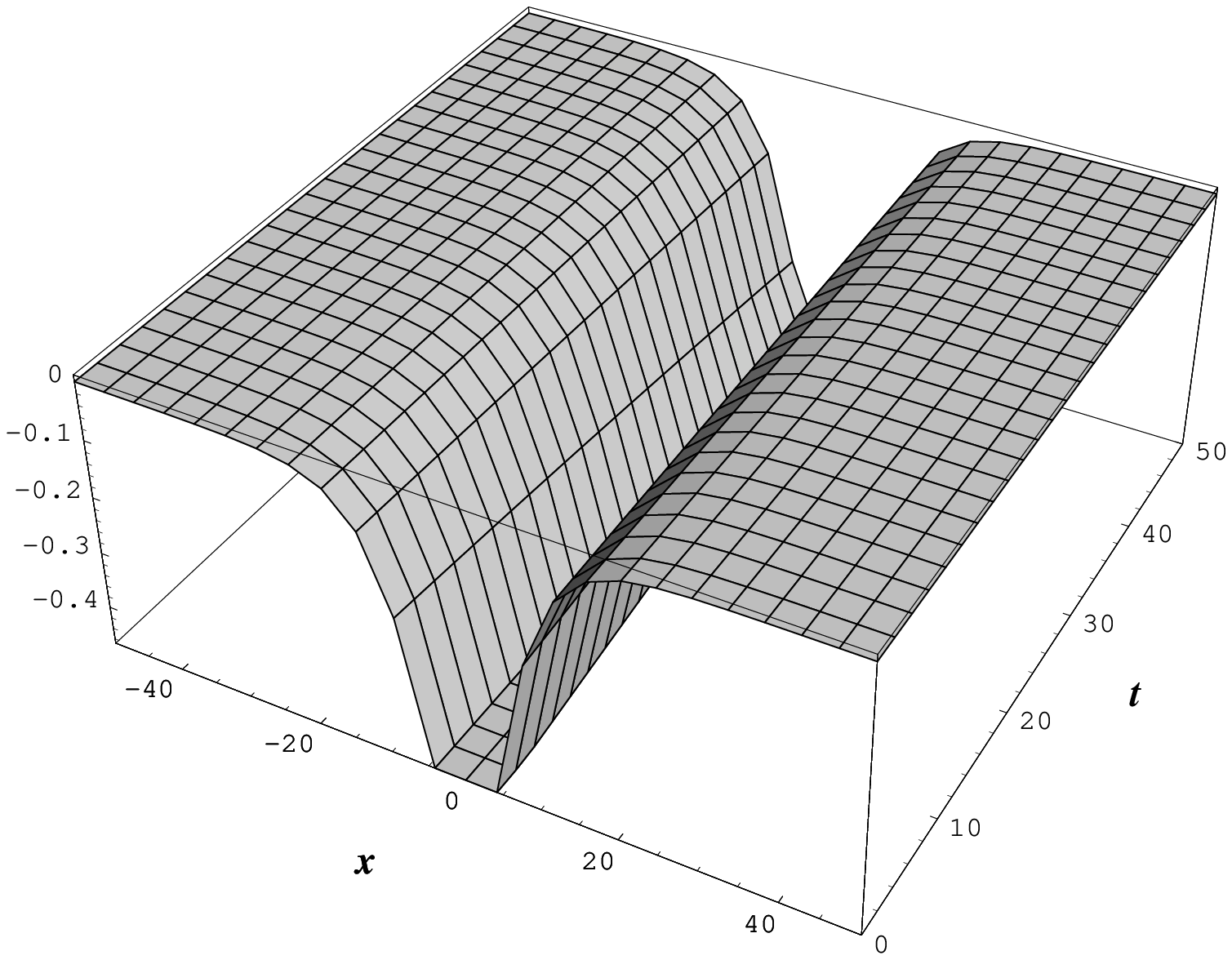}}
\caption{\small{The real and imaginary parts of the traveling wave
when $\alpha=0.05$, $\beta=-0.15$, $s=1$ and $\mu=0.5$}.}
\label{Profils}
\end{figure}

\noindent Figure \ref{Profils} presents the real and imaginary parts
of the traveling solitary wave as functions of the space variable
$x$ and the time variable $t$ for $\alpha=0.05$, $\beta=-0.15$,
$s=1$, $\mu=0.5$. For this set of values: \begin{equation}
\varepsilon_1=-\varepsilon_2=-\varepsilon_3=1
\end{equation}

\noindent The real part of the traveling wave presents a
kink-profile, while the imaginary part presents a bell-profile, as
it could be expected, since $|B_1|=|D_1|$, $D_1$ being purely
imaginary.
\\

Here we wish to clarify that the coefficients of the solutions in
[5, pp.428-429] are incorrect. By comparison with the previous
codes, we find that it is mainly due to Mathematica codes used for
[5], in which a cubic term was missed while authors typed the
algebraic system.

\section{Discussion}

\noindent In this work, by applying a hyperbolic ansatz, we obtain a
class of new kink-profile traveling solitary wave solutions to the
cBKdV equation (1), which indicates the coefficients of the
solutions in a previous work [5] are incorrect. One of the
advantages of this approach is that it is very straightforward to
generate new solutions and easily carried out with the aid of
mathematical softwares.
\\

\noindent Although explicit forms of traveling solitary wave
solutions are described herein and in the mentioned literature, as
far as our knowledge goes, it is still unclear whether the cBKdV
equation has any other type of bounded traveling solitary wave
solutions. In the near future, we are going to use Lie group
analysis as well as some other innovative methods to continue
studying this challenging problem. Physical explanation for the
solutions will be provided and the stability will be also analyzed
as the parameters vary.
\\

\addcontentsline{toc}{section}{\numberline{}References}


\begin{thebibliography}{1}

\normalsize{
\bibitem{burger1} Burgers J. M., Mathematical examples illustrating
relations occurring in the theory of turbulent fluid motion, {\em
Trans. Roy. Neth. Acad. Sci.} Amsterdam, 17 (1939) 1-53.

\bibitem{kdv1} Korteweg D. J. and de Vries G., On the change of
form of long waves advancing in a rectangular channel, and on a
new type of long stationary waves, {\em Phil. Mag.} 39 (1895)
422-443.

\bibitem{wadati} Wadati M., The modified Korteweg-de Vries equation,
{\em J. Phys. Soc. Japan}, 34 (1973) 1289-1296.

\bibitem{wang1} Wang M. L., Exact solutions for a compound KdV-Burgers equation,
{\em Phys. Lett. A,} 213 (1996) 279-287.

\bibitem{feng1} Feng Z. and Chen G., Solitary Wave Solutions of the Compound
Burgers-Korteweg-de Vries Equation, {\em Physica A}, 352 (2005)
419-435.

\bibitem{feng2} Feng, Z., A note on ``Explicit exact solutions to the compound
Burgers–Korteweg–de Vries equation", {\em Phys. Lett. A}, 312
(2003) 65-70.

\bibitem{feng3} Feng, Z., On explicit exact solutions to the compound Burgers-KdV equation,
{\em Phys. Lett. A}, 293 (2002) 57-66.

\bibitem{park1} Parkes E. J. and Duffy, B. R., Traveling solitary wave solutions to a compound KdV-Burgers
equation, {\em Phys. Lett. A} 229 (1997) 217-220.

\bibitem{park2} Parkes E. J., A note on solitary-wave solutions to compound KdV–Burgers
equations, {\em Phys. Lett. A} 317 (2003) 424-428.

\bibitem{zhang1} Zhang W. G., Chang Q. S. and Jiang B. G., Explicit exact solitary-wave solutions
for compound KdV-type and compound KdV–Burgers-type equations with
nonlinear terms of any order, {\em Chaos, Solitons \&  Fractals},
13 (2002) 311-319.

\bibitem{zhang2} Zhang W. G., Exact solutions of the Burgers–combined KdV mixed
equation, {\em Acta Math. Sci.} 16 (1996) 241–248.

\bibitem{li} Li B., Chen Y. and Zhang H. Q., Explicit exact solutions for new general
two-dimensional KdV-type and two-dimensional KdV–Burgers-type
equations with nonlinear terms of any order, {\em J. Phys. A
(Math. Gen.)} 35 (2002) 8253–8265.

\bibitem{whitham} Whitham G. B., Linear and Nonlinear Wave,
Wiley-Interscience, New York, 1974.

\bibitem{ablowitz} Ablowitz M. J. and Segur H., Solitons and the Inverse Scattering
Transform, SIAM, Philadelphia, 1981.


\bibitem{dodd} Dodd R. K., Eilbeck J. C., Gibbon J.D. and Morris H. C., Solitons and Nonlinear Wave
Equations, London Academic Press, London, 1983.

\bibitem{johnson} Johnson R. S., A Modern Introduction to the Mathematical Theory of
Water Waves, Cambridge University Press, Cambridge, 1997.

\bibitem{ince} Ince E.L., Ordinary Differential Equations, Dover
Publications, New York, 1956.

\bibitem{zhang3} Zhang Z. F., Ding T.R., Huang W. Z. and Dong Z. X.,
Qualitative Analysis of Nonlinear Differential Equations, Science
Press, Beijing, 1997.

%\bibitem{birk} Birkhoff G. and Rota G. C., Ordinary Differential
%Equations, Wiley, New York, 1989.
}


\end{thebibliography}
\end{document}